# R-Estimates vs. GMM: A Theoretical Case Study of Validity and Efficiency

**Dylan S. Small, Joseph L. Gastwirth, Abba M. Krieger and Paul R. Rosenbaum**

*Abstract.* What role should assumptions play in inference? We present a small theoretical case study of a simple, clean case, namely the nonparametric comparison of two continuous distributions using (essentially) information about quartiles, that is, the central information displayed in a pair of boxplots. In particular, we contrast a suggestion of John Tukey—that the validity of inferences should not depend on assumptions, but assumptions have a role in efficiency—with a competing suggestion that is an aspect of Hansen's generalized method of moments—that methods should achieve maximum asymptotic efficiency with fewer assumptions. In our case study, the practical performance of these two suggestions is strikingly different. An aspect of this comparison concerns the unification or separation of the tasks of estimation assuming a model and testing the fit of that model. We also look at a method (MERT) that aims not at best performance, but rather at achieving reasonable performance across a set of plausible models.

*Dylan Small is Assistant Professor, Department of Statistics, The Wharton School, University of Pennsylvania, Philadelphia, Pennsylvania 19104-6340, USA e-mail: dsmall@wharton.upenn.edu. Joseph L. Gastwirth is Professor of Statistics and Economics, George Washington University, Washington, District of Columbia 20052, USA. Abba M. Krieger is Robert Steinberg Professor of Statistics, Operations and Information Management, and Marketing, and Paul R. Rosenbaum is Robert G. Putzel Professor of Statistics, Department of Statistics, The Wharton School, University of Pennsylvania, Philadelphia, Pennsylvania 19104-6340, USA.*

## 1. INTRODUCTION: A QUESTION AND AN EXAMPLE

### 1.1 What Role for Assumptions?

In his essay, "Sunset Salvo," Tukey (1986, page 72) advocated:

> Reducing dependence on assumptions . . . using assumptions as leading cases, not truths, . . . when possible, using randomization to ensure validity—leaving to assumptions the task of helping with stringency.

Although the comment is not formal, presumably "validity" refers to the level of tests and the coverage rate of confidence intervals, while "stringency" refers to efficiency at least against some alternatives. Later in the essay (page 73), Tukey describes a statistic as "safe" if it is "valid—and of reasonably high efficiency—in each of a variety of situations." [Recall that a most stringent test minimizes the maximum power loss and that in many prob-lems, uniformly most powerful tests are not available; see Lehmann (1997).] The first part of Tukey's suggestion—"reduce dependence of validity on assumptions"—is today uncontroversial and there are many widely varied attempts to achieve that goal. The second part of Tukey's suggestion—"use assumptions to help with stringency"—runs against the grain of some recent developments which attempt to reduce the role of assumptions in obtaining efficient procedures. Do assumptions have a role in efficiency when comparing equally valid procedures? Or can we have it all, asymptotically of course, adapting our procedures to the data at hand to increase efficiency? Our purpose here is to closely examine these questions in a theoretical case study of a simple, clean case. We offer exactly the same information to two types of nonparametric procedures in a setting in which both are valid, though one chooses a procedure with high relative efficiency across a set of plausible models, while the other tries to be asymptotically efficient with fewer assumptions. The first method uses a form of rank statistic (Gastwirth, 1966, 1985; Birnbaum and Laska, 1967). The second method is a particular case of Hansen's (1982) generalized method of moments (GMM), widely used in econometrics. Both methods compare two distributions to estimate a shift. Both methods look at exactly the same information, somewhat related to the information about quartiles depicted in a pair of boxplots, but the methods use this simple information very differently. We compare the methods in a scientific example, in a simulation and using asymptotics. Also, we ask whether there is information against the shift model. We also show how to eliminate a shared assumption of both methods, namely the existence of a shift, which if false may invalidate their conclusions.

In Section 1.2 a motivating example is described. In Section 2 notation and the methods of estimating a shift are defined, and in Section 2.4 they are applied to the motivating example. The methods are evaluated by simulation in finite samples in Section 3, where some large sample results hold in quite small samples and others require astonishingly large samples. In Section 4 we dispense with the shift model. The relevant large sample theory is discussed in the Appendix with some patches needed to cover some nonstandard details.

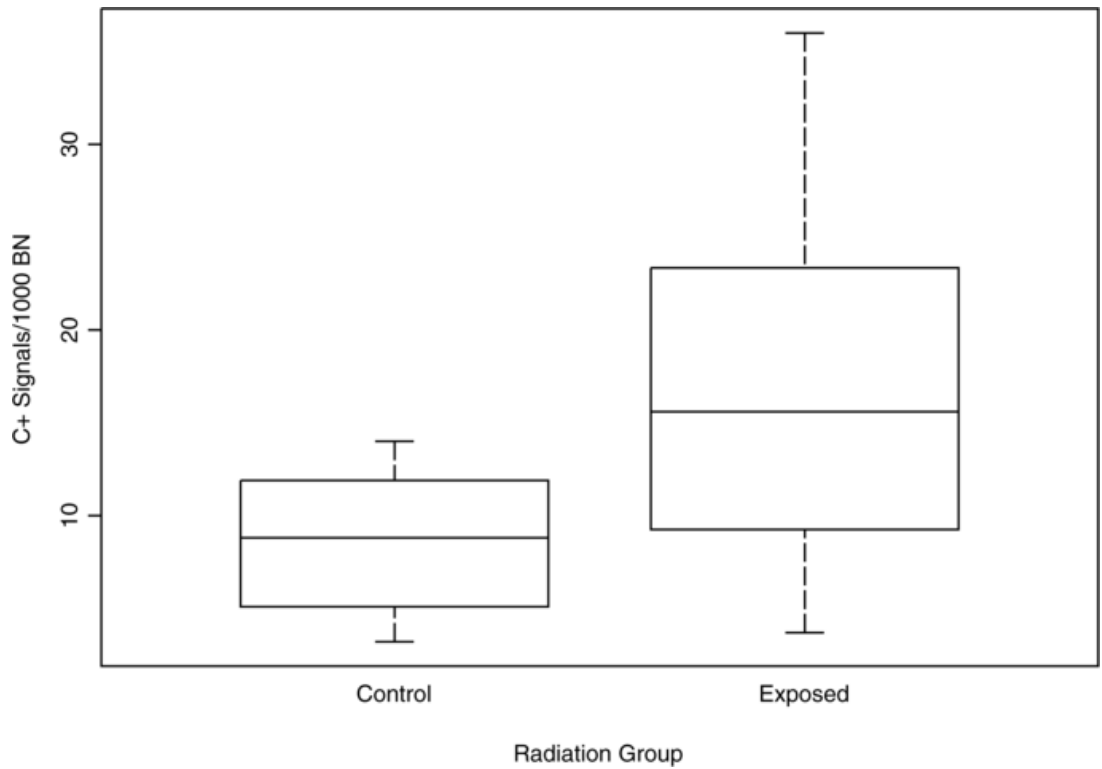

Fig. 1. *Genetic damage in radiation exposed and control groups.*

## 1.2 A Motivating Example: Radiation in Homes

In the early 1980s, a number of residential buildings were constructed in Taiwan using $^{60}Co$-contaminated steel rods, with the consequence that the levels of radioactive exposure in these homes were often orders of magnitude higher than background levels. Chang et al. (1999) compared 16 residents of these buildings to 7 unexposed controls with respect to several measures of genetic damage, including the number of centromere-positive signals per 1000 binucleated cells, as depicted in Figure 1. The sorted values for the 16 residents were 3.7, 6.8, 8.4, 8.5, 10.0, 11.3, 12.0, 12.5, 18.7, 19.0, 20.0, 22.7, 24.0, 31.8, 33.3, 36.0 and for the 7 controls were 3.2, 5.1, 8.3, 8.8, 9.5, 11.9, 14.0. They reported means, standard deviations and the significance level from the Mann–Whitney–Wilcoxon test. In Figure 1 the distribution for exposed subjects looks higher, more dispersed and possibly slightly skewed right in comparison to the controls, but of course the sample sizes are small and boxplots fluctuate in appearance by chance. Should one estimate a shift in the distributions? If so, how? If not, what should one do instead?

Suppose that the control measurements, $X_1, \ldots, X_m$, are independent observations from a continuous, strictly increasing cumulative distribution $F(\cdot)$ and that the exposed measurements, $Y_1, \ldots, Y_n$, are independent observations from a continuous, strictly increasing cumulative distribution $G(\cdot)$, with $N = n + m$. The distributions are shifted if there is some constant $\Delta$ such that $Y_j - \Delta$ and $X_i$ have the same distribution $F(\cdot)$, or equivalently if $F(x) = G(x+\Delta)$ for each $x$. We are interested in whether a shift model is compatible with the data, and if it is, in estimating $\Delta$, and if it is not, in estimating something else. Obviously, in the end, with finite amounts of data, there is going to be some uncertainty about both questions—whether a shift model is appropriate, what values of $\Delta$ are reasonable—but the goal is to describe the available information and the uncertainty that remains.

## 2. TWO APPROACHES USING THE SAME INFORMATION

### 2.1 A $4 \times 2$ Table

Boxplots serve two purposes: they call attention to unique, extreme observations requiring individual attention, and they describe distributional shape in terms of quartiles. Here we focus on quartiles, and contrast two approaches to using (more or less) the depicted information to determine whether a shift, $\Delta$, exists and what values of $\Delta$ are reasonable. We very much want to offer the two approaches exactly the same information, and then see which approach makes better use of this information; that is, we want to avoid complicating the comparison by offering different information to the different approaches. Because the question about distributions is raised by the appearance of boxplots, the information offered to the methods is information about the quartiles.

Consider the null hypothesis that the distributions are shifted by a specified amount, $\Delta_0$, that is, $H_0 : F(x) = G(x + \Delta_0)$. If this hypothesis were true, then $Z_1^{\Delta_0} = X_1, \ldots, Z_m^{\Delta_0} = X_m, Z_{m+1}^{\Delta_0} = Y_1 - \Delta_0, \ldots, Z_N^{\Delta_0} = Y_n - \Delta_0$, would be $N$ independent observations from $F(\cdot)$. Let $Z_{(1)}^{\Delta_0} < \cdots < Z_{(N)}^{\Delta_0}$ be the order statistics, let $q_i = \lceil \frac{iN}{4} \rceil$ for $i = 1, 2, 3$, where $\lceil w \rceil$ is the least integer greater than or equal to $w$, and define $Z_{(q_i)}^{\Delta_0}$ to be the $i$th quartile. With $N = 23$ in the example, $q_1 = 6$, $q_2 = 12$, $q_3 = 18$ and the quartiles are $Z_{(6)}^{\Delta_0}$, $Z_{(12)}^{\Delta_0}$, $Z_{(18)}^{\Delta_0}$. Write $k_1 = q_1$, $k_2 = q_2 - q_1$, $k_3 = q_3 - q_2$ and $k_4 = N - q_3$, so for $N = 23$, $k_1 = 6$, $k_2 = 6$, $k_3 = 6$, $k_4 = 5$. For the hypothesized $\Delta_0$, form the $4 \times 2$ contingency table in Table 1 which classifies the $Z_i^{\Delta_0}$ by quartile and by treatment or control. Notice that the marginal totals of Table 1 are functions of the sample sizes, $n$, $m$ and $N = n + m$, so they are fixed, not varying from sample to sample. If the null hypothesis is true, Table 1 has the multivariate hypergeometric distribution

$$
\begin{aligned}
&\Pr(A_1^{\Delta_0} = a_1, A_2^{\Delta_0} = a_2, A_3^{\Delta_0} = a_3, A_4^{\Delta_0} = a_4) \\
&= \frac{\binom{k_1}{a_1}\binom{k_2}{a_2}\binom{k_3}{a_3}\binom{k_4}{a_4}}{\binom{N}{n}}
\end{aligned}
$$
(1)

for $0 \le a_j \le k_j$, $j = 1, 2, 3, 4$, $a_1 + a_2 + a_3 + a_4 = n$.

TABLE 1
*Contingency table from pooled quartiles*

| Quartile interval | Treated | Control | Total |
|---|---|---|---|
| $Z_i^{\Delta_0} \le Z_{(q_1)}^{\Delta_0}$ | $A_1^{\Delta_0}$ | $k_1 - A_1^{\Delta_0}$ | $k_1$ |
| $Z_{(q_1)}^{\Delta_0} < Z_i^{\Delta_0} \le Z_{(q_2)}^{\Delta_0}$ | $A_2^{\Delta_0}$ | $k_2 - A_2^{\Delta_0}$ | $k_2$ |
| $Z_{(q_2)}^{\Delta_0} < Z_i^{\Delta_0} \le Z_{(q_3)}^{\Delta_0}$ | $A_3^{\Delta_0}$ | $k_3 - A_3^{\Delta_0}$ | $k_3$ |
| $Z_{(q_3)}^{\Delta_0} < Z_i^{\Delta_0}$ | $A_4^{\Delta_0}$ | $k_4 - A_4^{\Delta_0}$ | $k_4$ |
| Total | $n$ | $m$ | $N$ |

As a consequence, if the null hypothesis is true, the expected counts are

$$E(A_j^{\Delta_0}) = \frac{nk_j}{N}, \quad j = 1, 2, 3, 4, \tag{2}$$

with variances and covariances

$$\operatorname{var}(A_j^{\Delta_0}) = \frac{nmk_j(N - k_j)}{N^2(N - 1)},$$
$$\operatorname{cov}(A_i^{\Delta_0}, A_j^{\Delta_0}) = -\frac{nmk_ik_j}{N^2(N - 1)}, \quad i \neq j. \tag{3}$$

Write

$$\mathbf{A}_{\Delta_0} = (A_1^{\Delta_0}, A_2^{\Delta_0}, A_3^{\Delta_0}, A_4^{\Delta_0})^T,$$
$$\mathbf{E} = \left(\frac{nk_1}{N}, \frac{nk_2}{N}, \frac{nk_3}{N}, \frac{nk_4}{N}\right)^T$$

and $\mathbf{V}$ for the symmetric $4 \times 4$ covariance matrix of the hypergeometric defined by (3). Notice, in particular, that $\mathbf{A}_{\Delta_0}$ is a random vector whose value changes with $\Delta_0$, whereas $\mathbf{E}$ and $\mathbf{V}$ are fixed matrices whose values do not change with $\Delta_0$. For each value of $\Delta_0$ there is a table of the form Table 1, and if the distributions were actually shifted, $F(x) = G(x + \Delta)$, then the table with $\Delta_0 = \Delta$ has the multivariate hypergeometric distribution. The information available to both methods of inference is this collection of $4 \times 2$ tables as $\Delta_0$ varies.

A minor technical issue needs to be mentioned. Because $A_1^{\Delta_0} + A_2^{\Delta_0} + A_3^{\Delta_0} + A_4^{\Delta_0} = n$ is constant, the covariance matrix $\mathbf{V}$ is singular, that is, positive semidefinite but not positive definite. One could avoid this by focusing on $(A_2^{\Delta_0}, A_3^{\Delta_0}, A_4^{\Delta_0})$, but the $4 \times 2$ table is too familiar to discard for this minor technicality. We define $\mathbf{V}^-$ as the specific generalized inverse of $\mathbf{V}$ which has 0's in its first row and column, and has in its bottom right $3 \times 3$ corner the inverse of the bottom right $3 \times 3$ corner of $\mathbf{V}$; see Rao (1973, page 27). Although our notation always refers to the $4 \times 2$ table, the calculations ultimately use only the nondegenerate piece of the table, $(A_2^{\Delta_0}, A_3^{\Delta_0}, A_4^{\Delta_0})$. In later sections, this issue comes up several times, always with minor consequences.

The distribution (1) for Table 1 also arises in other ways. For instance, if there are $N$ subjects and $n$ are randomly assigned to treatment, with the remaining $m = N - n$ subjects assigned to control, and if the treatment has an additive effect $\Delta_0$, then (1) is the distribution for Table 1 without assuming samples from a population.

## 2.2 Inference Based on Group Ranks

A simple approach to inference about $\Delta$ assuming shifted distributions uses a group rank statistic, as discussed by Gastwirth (1966) and Markowski and Hettmansperger (1982, Section 5); see also Brown (1981) and Rosenbaum (1999). Here, the rows of Table 1 are assigned scores, $\mathbf{w} = (w_1, w_2, w_3, w_4)^T$, with $w_1 = 0$, and the hypothesis $H_0: \Delta = \Delta_0$ is tested using the group rank statistic $T_{\Delta_0} = \mathbf{w}^T \mathbf{A}_{\Delta_0}$ whose exact null distribution is determined using (1). An exact confidence interval for $\Delta$ is obtained by inverting the test; for example, see Lehmann (1963), Moses (1965) and Bauer (1972). The Hodges–Lehmann (1963) point estimate, $\widehat{\Delta}_{\mathrm{HL}}$, for this rank test is essentially the solution to the estimating equation

$$\mathbf{w}^T \mathbf{A}_{\widetilde{\Delta}} = \mathbf{w}^T \mathbf{E}. \tag{4}$$

More precisely, with increasing scores, $0 = w_1 < w_2 < w_3 < w_4$, because $\mathbf{w}^T \mathbf{A}_{\widetilde{\Delta}}$ moves in discrete steps as $\widetilde{\Delta}$ varies continuously, the Hodges–Lehmann estimate is defined so that: (i) if equality in (4) cannot be achieved, then the estimate is the unique point $\widetilde{\Delta}$ where $\mathbf{w}^T \mathbf{A}_{\widetilde{\Delta}}$ passes $\mathbf{w}^T \mathbf{E}$, or (ii) if equality is achieved for an interval of values of $\widetilde{\Delta}$ then $\widehat{\Delta}_{\mathrm{HL}}$ is defined to be the midpoint of the interval. In large samples, the null distribution of $T_{\Delta_0}$ is approximately Normal with expectation $\mathbf{w}^T \mathbf{E}$ and variance $\mathbf{w}^T \mathbf{V} \mathbf{w}$, so the deviate

$$D_{\Delta_0} = \frac{\mathbf{w}^T(\mathbf{A}_{\Delta_0} - \mathbf{E})}{\sqrt{\mathbf{w}^T \mathbf{V} \mathbf{w}}} \tag{5}$$

is compared with the standard Normal distribution. The Hodges–Lehmann estimate $\widehat{\Delta}_{\mathrm{HL}}$ defined by (4) is essentially the same as the value $\widetilde{\Delta}$ that minimizes $D_{\widetilde{\Delta}}^2$.

All of this assumes the distributions are indeed shifted. A simple test of the hypothesis that $F(x) = G(x + \Delta_0)$ is based on the statistic

$$G_{\Delta_0}^2 = (\mathbf{A}_{\Delta_0} - \mathbf{E})^T \mathbf{V}^-(\mathbf{A}_{\Delta_0} - \mathbf{E}),$$

whose exact null distribution follows from (1) and whose large sample null distribution is approximately $\chi^2$ on three degrees of freedom. It is useful to notice that under the null hypothesis $E(G_{\Delta_0}^2) = \operatorname{tr}[E\{(\mathbf{A}_{\Delta_0} - \mathbf{E})(\mathbf{A}_{\Delta_0} - \mathbf{E})^T\}\mathbf{V}^-] = 3$, so the exact null distribution of $G_{\Delta_0}^2$ and the $\chi^2$ approximation have the same expectation. This will be relevant to certain comparisons to be made later. Is the shift model plausible for plausible values of $\Delta_0$? A simple, informative procedure is to plot the exact, two-sided $P$-values from $D^2_{\Delta_0}$ and $G^2_{\Delta_0}$ against $\Delta_0$. Curiosity is, of course, aroused by values $\Delta_0$ which are accepted by $D^2_{\Delta_0}$ and rejected by $G^2_{\Delta_0}$, because then the shift model is implausible for an ostensibly plausible shift. Greevy et al. (2004) do something similar.

## 2.3 Inference Based on the Generalized Method of Moments

In an important and influential paper, Hansen (1982) proposed a method for combining a number of estimating equations to estimate a smaller number of parameters. In the current context, the Tables 1 for different $\Delta_0$ yield the four estimating equations given by (2), one of which is redundant or linearly dependent on the other three. Obviously, there may be no $\Delta_0$ that satisfies all four equations at once, so Hansen proposed weighting the equations in an optimal way. In particular, he showed that the optimal weighting of the equations uses the inverse covariance matrix of the moment conditions (2), and it is, in fact, the value $\widehat{\Delta}_{\mathrm{gmm}}$ that minimizes $G^2_{\Delta_0}$:

$$(6) \quad \widehat{\Delta}_{\mathrm{gmm}} = \underset{\widetilde{\Delta}}{\arg\min} (\mathbf{A}_{\widetilde{\Delta}} - \mathbf{E})^T \mathbf{V}^- (\mathbf{A}_{\widetilde{\Delta}} - \mathbf{E}).$$

In theory, $\widehat{\Delta}_{\mathrm{gmm}}$ is asymptotically efficient, fully utilizing all of the information in the estimating equations (2); see Hansen (1982) and the Appendix. Surveys of GMM are given by Mátyás (1999) and Lindsay and Qu (2003). Hansen's results do not quite apply here, because certain differentiability assumptions he makes are not strictly satisfied, but his conclusions hold nonetheless, as discussed in the Appendix where a result of Jurečková (1969) about the asymptotic linearity of rank statistics replaces differentiability. Hansen proposed a large sample test of the model or "identifying restrictions" (2) using the minimum value of $G^2_{\Delta_0}$, that is, here, testing the family of models by comparing

$$G^2_{\widehat{\Delta}_{\mathrm{gmm}}} = (\mathbf{A}_{\widehat{\Delta}_{\mathrm{gmm}}} - \mathbf{E})^T \mathbf{V}^- (\mathbf{A}_{\widehat{\Delta}_{\mathrm{gmm}}} - \mathbf{E})$$

to the $\chi^2$ distribution on two degrees of freedom. There does not appear to be an exact null distribution for $G^2_{\widehat{\Delta}_{\mathrm{gmm}}}$ because it is computed not at $\Delta$ but rather at $\widehat{\Delta}_{\mathrm{gmm}}$, which varies from sample to sample, so the hypergeometric distribution for $\mathbf{A}_\Delta$ is not relevant. One large sample test of $H_0 : \Delta = \Delta_0$ compares $G^2_{\Delta_0} - G^2_{\widehat{\Delta}_{\mathrm{gmm}}}$ to the chi-square distribution on one degree of freedom (Newey and West, 1987; Mátyás, 1999, page 109), and a confidence set is the set of $\Delta_0$ not rejected by the test. See the Appendix.

By a familiar fact (Rao, 1973, page 60)

$$(7) \quad \sup_{\mathbf{c} = (0, c_2, c_3, c_4)} \frac{\{\mathbf{c}^T (\mathbf{A}_{\widetilde{\Delta}} - \mathbf{E})\}^2}{\mathbf{c}^T \mathbf{V} \mathbf{c}}$$
$$= (\mathbf{A}_{\widetilde{\Delta}} - \mathbf{E})^T \mathbf{V}^- (\mathbf{A}_{\widetilde{\Delta}} - \mathbf{E}) = G^2_{\widetilde{\Delta}},$$

so that

$$(8) \quad \widehat{\Delta}_{\mathrm{gmm}} = \underset{\widetilde{\Delta}}{\arg\min} \sup_{\mathbf{c} = (0, c_2, c_3, c_4)} \frac{\{\mathbf{c}^T (\mathbf{A}_{\widetilde{\Delta}} - \mathbf{E})\}^2}{\mathbf{c}^T \mathbf{V} \mathbf{c}},$$

whereas $\widehat{\Delta}_{\mathrm{HL}}$ minimizes an analogous quantity, namely $D^2_{\widetilde{\Delta}}$ in (5) for one specific set of weights $\mathbf{w}$. The $\mathbf{w}$ that achieves the bound in (7) is $\mathbf{w} = \mathbf{V}^- (\mathbf{A}_{\widetilde{\Delta}} - \mathbf{E})$, so this optimizing $\mathbf{w}$ is not fixed, but is rather a function of the data.

Both $\widehat{\Delta}_{\mathrm{HL}}$ and $\widehat{\Delta}_{\mathrm{gmm}}$ use the same information, the information in Table 1 for varied $\Delta_0$ and the moment equations (2); however, $\widehat{\Delta}_{\mathrm{HL}}$ uses an a priori weighting of the equations yielding the one estimating equation (4), while $\widehat{\Delta}_{\mathrm{gmm}}$ weights the equations (2) using $\mathbf{V}^-$. In the sense of (7), $\widehat{\Delta}_{\mathrm{gmm}}$ uses the "best" weights as judged by the sample, and asymptotically it achieves the efficiency associated with knowing what are the best fixed weights to use; see the Appendix. Is it best to use the "best" weights?

The generalized method of moments includes many familiar methods of estimation, including maximum likelihood, least squares and two-stage least squares with instrumental variables. In particular cases, such as weak instruments, it is known that poor estimates may result from GMM (e.g., Imbens, 1997; Staiger and Stock, 1997), but this is sometimes viewed as a weakness in the available data. A weak instrument may reduce efficiency but need not result in invalidity if appropriate methods of analysis are used (Imbens and Rosenbaum, 2005). The two-sample shift problem is identified and presents no weakness in the data.

## 2.4 Example

The methods of Sections 2.2 and 2.3 will now be applied to the data in Figure 1, where the two boxplots have medians 8.8 and 15.6 differing by 6.8, and means 8.7 and 17.4 differing by 8.7. If the distributions were shifted by $\Delta$, then $T_\Delta$ would be distribution free using (1), and the expectation of $T_\Delta$ with rank weights, $w_j = j$, $j = 0, 1, 2, 3$, would be 22.957 and the variance would be 6.1206. Now, $T_{8.7} = 22$, $T_{8.69999} = 23$, so $\widehat{\Delta}_{\mathrm{HL}} = 8.7$, which is by coincidence the same as the difference in means. Also, $G^2_{\Delta_0}$ takes

## Shifted Boxplots Using HL and GMM

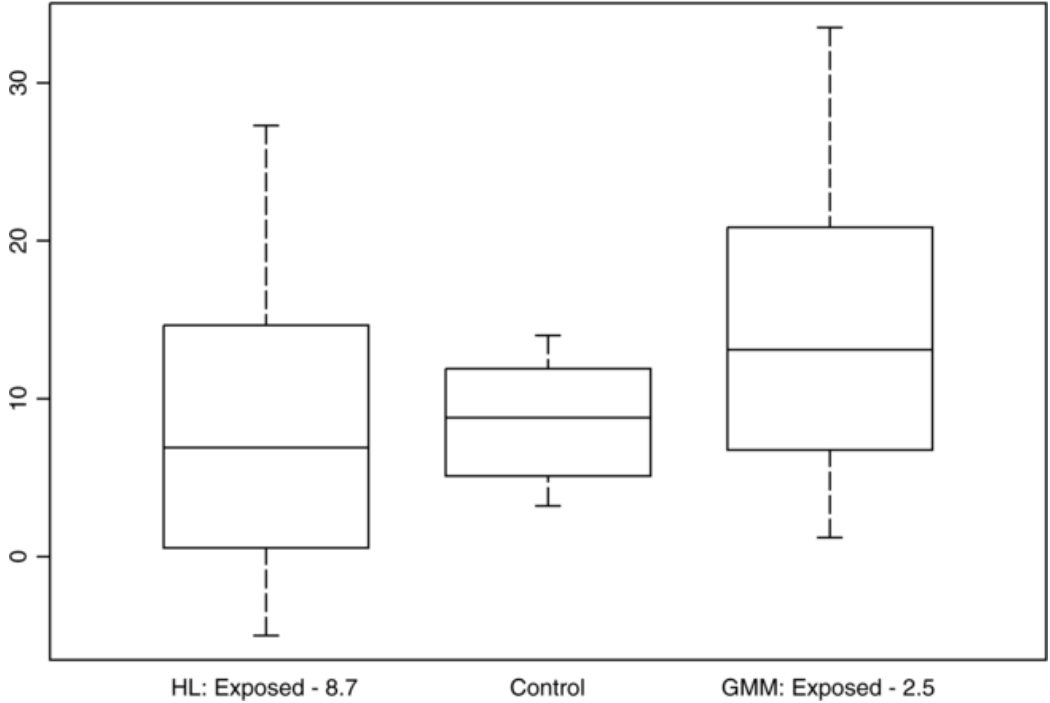

FIG. 2.   *Boxplots of $Y_j - \widehat{\Delta}_{\mathrm{HL}}$, $X_i$ and $Y_j - \widehat{\Delta}_{\mathrm{gmm}}$. GMM failed to align the boxes: all of the quartiles of $Y_j - \widehat{\Delta}_{\mathrm{gmm}}$ are above those of $X_i$.*

its smallest value, 3.176, on the half open interval $[0.10, 4.90)$, so $\widehat{\Delta}_{\mathrm{gmm}} = \frac{0.10 + 4.90}{2} = 2.5$. In these data, which estimate, $\widehat{\Delta}_{\mathrm{HL}} = 8.7$ or $\widehat{\Delta}_{\mathrm{gmm}} = 2.5$, looks better?

Figure 2 compares $\widehat{\Delta}_{\mathrm{HL}}$ and $\widehat{\Delta}_{\mathrm{gmm}}$ by displaying the control responses, $X_i$, together with the adjusted exposed responses, $Y_j - \widehat{\Delta}_{\mathrm{HL}}$ or $Y_j - \widehat{\Delta}_{\mathrm{gmm}}$. If the shift model were true, the estimates equaled the true shift, and the sample size were very large, then the three boxplots would look essentially the same. As it is, $Y_j - \widehat{\Delta}_{\mathrm{gmm}}$ appears both too high and too dispersed compared to the controls: all three quartiles of the $Y_j - \widehat{\Delta}_{\mathrm{gmm}}$ are above the corresponding quartiles of the $X_i$; the median of $Y_j - \widehat{\Delta}_{\mathrm{gmm}}$ is above the upper quartile of the $X_i$; the upper quartile of $Y_j - \widehat{\Delta}_{\mathrm{gmm}}$ is above the maximum of $X_i$. In contrast, the $Y_j - \widehat{\Delta}_{\mathrm{HL}}$ are shifted reasonably but appear more dispersed than the $X_i$: the upper quartile of the $Y_j - \widehat{\Delta}_{\mathrm{HL}}$ is too high while the lower quartile is too low; indeed, the entire boxplot of the $X_i$ fits inside the quartile box of the $Y_j - \widehat{\Delta}_{\mathrm{HL}}$. Obviously, a shift can relocate a boxplot but cannot alter its dispersion.

Assuming there is a shift $\Delta$, the exact distribution of the squared deviate $D_\Delta^2$ based on $T_\Delta$ is determined using (1) and it has $\Pr(D_\Delta^2 \geq 4.156) = 0.0436$. The $1 - 0.0436 = 95.6\%$ confidence set $\mathcal{D}$ is the closure of the set $\{\Delta_0 : D_{\Delta_0}^2 < 4.156\}$, which is $[0.1, 19.5]$. The test of fit of the shift model based on the generalized method of moments yields $G_{\widehat{\Delta}_{\mathrm{gmm}}}^2 = G_{2.5}^2 = 3.176$, which is compared to $\chi^2$ on two degrees of freedom, yielding a significance level greater than 0.2. In short, the GMM test of the shift model based on $G_{\widehat{\Delta}_{\mathrm{gmm}}}^2$ suggests the shift model is plausible, and the appearance of Figure 2 could be due to chance. Figure 3 is the plot, suggested in Section 2.2, of the exact, two-sided $P$-values from $D_{\Delta_0}^2$ and $G_{\Delta_0}^2$ plotted against $\Delta_0$. To focus attention on small $P$-values, the vertical axis uses a logarithmic scale. To anchor that scale, horizontal lines are drawn at $P = 0.05$, 0.1 and 1/3. The solid step function for $D_{\Delta_0}^2$ cuts the horizontal $P = 0.05$ line at the endpoints for the 95% confidence interval. As suggested by Mosteller and Tukey (1977), a 2/3 confidence interval is analogous to an estimate plus or minus a standard error. The Hodges–Lehmann estimate is $\widehat{\Delta}_{HL} = 8.7$. Notice, however, that at $\Delta_0 = 8.69$, the $P$-value for $D^2_{\Delta_0}$ is 1.000, but the $P$-value for $G^2_{\Delta_0}$ is 0.021; that is, the shift model is implausible for a value of $\Delta_0$ judged highly plausible by $D^2_{\Delta_0}$. Table 2 is Table 1 evaluated at $\Delta_0 = 8.69$, and from this table, it is easy to see what has happened. With $\Delta_0 = 8.69$ subtracted from treated responses, $T_{8.69} = 5 \times 3 + 3 \times 2 + 2 \times 1 + 0 \times 6 = 23$, which is as close as possible to the null expectation 22.957 of $T_\Delta$, but because of the greater dispersion in the treated group, all $6 + 5 = 11$ observations outside the pooled upper and lower quartiles are treated responses, leading to a large $G^2_{8.69} = 9.2$ with exact significance level 0.021. The pattern in Table 2 is hardly a surprise: the comparison of dispersions is most decisive when it is not obscured by unequal locations.

Because $T_{\Delta_0}$ is monotone in $\Delta_0$, the 95%, 90% and 2/3 confidence sets it yields in Figure 3 are intervals. In contrast, the 95%, 90% and 2/3 confidence sets based on $G^2_{\Delta_0}$ are not intervals; for instance, the 90% confidence set is the union of three disjoint intervals. If the confidence interval is defined as the shortest closed interval containing the confidence set, then the three intervals based on $G^2_{\Delta_0}$ are all longer than the corresponding intervals based on $D^2_{\Delta_0}$. Figure 4 calculates the large sample confidence interval from GMM, plotting $G^2_{\Delta_0} - G^2_{\widehat{\Delta}_{gmm}}$ against $\Delta_0$. For instance, the dotted line labeled 95% in Figure 4 is at 3.841, the 95% point of the chi-square distribution with one degree of freedom. The 95% confidence set for $\Delta$ is the set of $\Delta_0$ such that $G^2_{\Delta_0} - G^2_{\widehat{\Delta}_{gmm}} \leq 3.841$, and it is the union of two disjoint intervals; similarly, the 90% confidence set is the union of three disjoint intervals, and the 2/3 confidence set is the union of two disjoint intervals. The shortest closed interval containing the 95% confidence set for $\Delta$ is $[-2.7, 19.5]$, which is longer than the exact 95% confidence interval based on $D^2_{\Delta_0}$, namely $[0.1, 19.5]$. Of course, both confidence intervals include values of $\Delta$ rejected by $G^2_{\Delta_0}$, so the shift model is not really plausible for some parameter values that both intervals report as plausible shifts.

In short, the method of Section 2.2 gave a point estimate of $\widehat{\Delta}_{HL} = 8.7$ consistent with the difference in means, but raised doubts about whether a shift model is appropriate, rejecting the shift model for $H_0 : \Delta = 8.69$. In contrast, the method of Section 2.3 suggested the shift is much smaller, $\widehat{\Delta}_{gmm} = 2.5$, and the associated goodness of fit test based on $G^2_{\widehat{\Delta}_{gmm}}$ suggested that the shift model is plausible. Obviously, the example just illustrates what the two methods do with one data set; it tells us nothing about performance in large or small samples.

## 3. SIMULATION

### 3.1 Structure of the Simulation

The simulation considered three distributions $F(\cdot)$, namely the Normal (N), the Cauchy (C) and the convolution of a standard Normal with a standard Exponential (NE). Recall that the standard Normal and Exponential distributions each have variance one, so NE has variance two. Although the support of NE is the entire line, NE has a long right tail and a short left tail, and is moderately asymmetric near its median. There were 5,000 samples drawn for each sampling situation.

We considered several estimators, including the two in Section 2.4, namely $\widehat{\Delta}_{gmm}$ based on GMM and $\widehat{\Delta}_{HL}$ using scores $w_1 = 0$, $w_2 = 1$, $w_3 = 2$, $w_4 = 3$. Note that the weights for $\widehat{\Delta}_{HL}$ are close to the optimal weights for the Normal. The estimate $\widehat{\Delta}_M$ with scores $w_1 = 0$, $w_2 = 0$, $w_3 = 1$, $w_4 = 1$ is the Hodges–Lehmann point estimate associated with Mood's two sample median test; these scores are close to the optimal scores for the Cauchy. The estimate $\widehat{\Delta}_{mert}$ with weights $w_1 = 0$, $w_2 = 0.18$, $w_3 = 0.82$, $w_4 = 1$ is Gastwirth's compromise weights for the Normal and the Cauchy; these scores are almost the same as $w_1 = 0$, $w_2 = 1$, $w_3 = 4$, $w_4 = 5$, so are much closer to $\widehat{\Delta}_M$ than to $\widehat{\Delta}_{HL}$. The coverage rates and behavior of confidence intervals and the null distribution of the goodness of fit test based on $G^2_{\widehat{\Delta}_{gmm}}$ were also examined.

### 3.2 Efficiency

Asymptotic theory says: (i) $\widehat{\Delta}_{gmm}$ should always win in sufficiently large samples, (ii) $\widehat{\Delta}_{HL}$ should be close to the best for the Normal, (iii) $\widehat{\Delta}_M$ should be close to the best for the Cauchy and (iv) $\widehat{\Delta}_{mert}$ should be better than $\widehat{\Delta}_{HL}$ for the Cauchy and better than $\widehat{\Delta}_M$ for the Normal.

TABLE 2
*Table for testing* 8.69

| Quartile interval | Treated | Control | Total |
|---|---|---|---|
| Lowest | 6 | 0 | 6 |
| Low | 2 | 4 | 6 |
| High | 3 | 3 | 6 |
| Highest | 5 | 0 | 5 |
| Total | 16 | 7 | 23 |

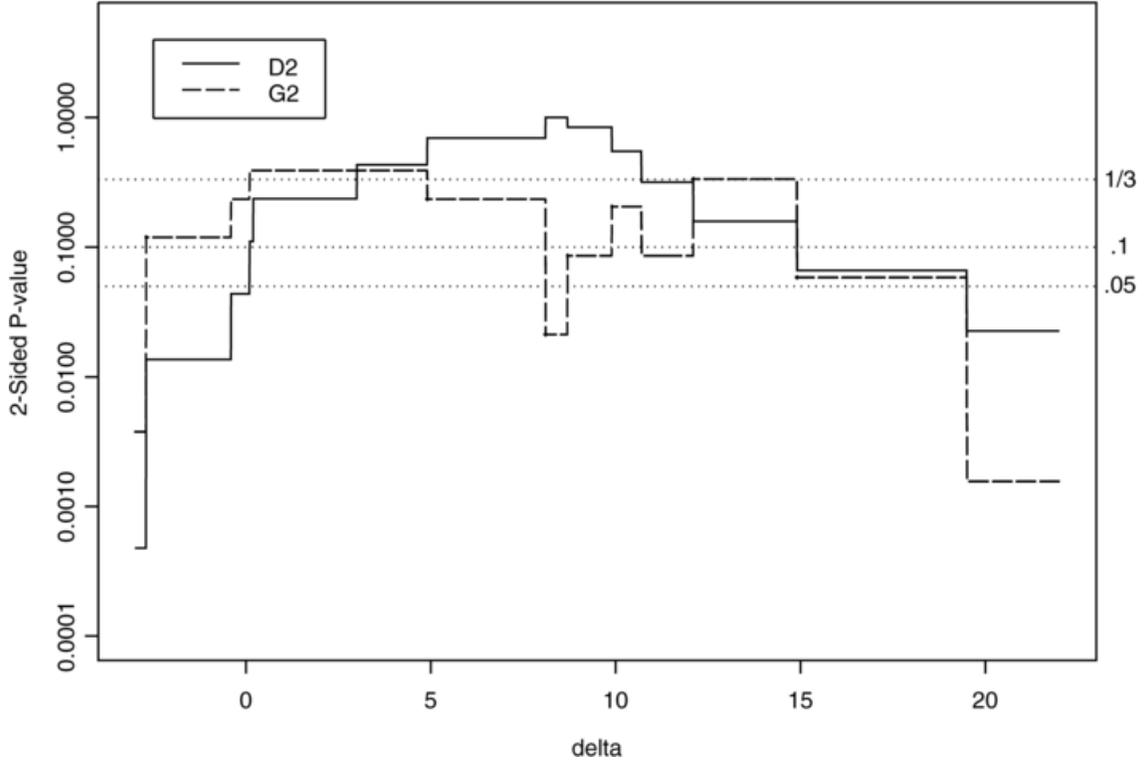

FIG. 3. *Plot of exact significance levels for testing $H_0 : \Delta = \Delta_0$ using $D^2_{\Delta_0}$ and $G^2_{\Delta_0}$ versus $\Delta_0$. Note that $G^2_{\Delta_0}$ rejects the $\Delta_0$ that minimizes $D^2_{\Delta_0}$.*

Table 3 compares efficiency for samples from the Normal distribution. The values in the table are ratios of mean squared errors averaged over 5000 samples, so the value 0.72 for $n = m = 24$ in $\widehat{\Delta}_{gmm} : \widehat{\Delta}_{HL}$ indicates that $\widehat{\Delta}_{HL}$ had a mean squared error that was 72% of the mean squared error of $\widehat{\Delta}_{gmm}$. The

TABLE 3
*Efficiency for samples from the Normal distribution*

| $n$ | 24 | 50 | 20 | 80 | 500 | 2,000 | 10,000 |
| $m$ | 24 | 50 | 80 | 80 | 500 | 2,000 | 10,000 |
|---|---|---|---|---|---|---|---|
| $\widehat{\Delta}_{gmm} : \widehat{\Delta}_{HL}$ | 0.72 | 0.76 | 0.81 | 0.79 | 0.86 | 0.89 | 0.93 |
| $\widehat{\Delta}_{gmm} : \widehat{\Delta}_M$ | 1.00 | 1.05 | 1.08 | 1.07 | 1.14 | 1.18 | 1.24 |
| $\widehat{\Delta}_{gmm} : \widehat{\Delta}_{mert}$ | 0.80 | 0.86 | 0.91 | 0.89 | 0.94 | 0.97 | 1.03 |

Summary: $\widehat{\Delta}_{HL}$ is best in all cases, and $\widehat{\Delta}_{mert}$ is second for $n, m \leq 2000$.

predictions of large sample theory are qualitatively correct, but some of the quantitative results are striking. The performance of $\widehat{\Delta}_{gmm}$ improves with increasing sample size, but $\widehat{\Delta}_{gmm}$ is still 7% behind $\widehat{\Delta}_{HL}$ with $n = m = 10,000$ and only marginally better than $\widehat{\Delta}_{mert}$. For sample sizes of $n = m = 2,000$ or less, $\widehat{\Delta}_{mert}$ is better than $\widehat{\Delta}_{gmm}$, often substantially so. In Table 3, the relative performance of $\widehat{\Delta}_{gmm}$ is still improving as the sample size increases, with the promised optimal performance not yet visible for the sample sizes in the table. With $n = m = 40,000$, not shown in the table, the relative efficiency $\widehat{\Delta}_{gmm}$ is still about 5% behind $\widehat{\Delta}_{HL}$.

Table 4 is the analogous table for samples from the Cauchy distribution. As before, the relative performance of $\widehat{\Delta}_{gmm}$ improves with increasing sample size, so that it is inferior to $\widehat{\Delta}_{HL}$ for $n = m = 80$ but superior for $n = m = 500$. In Table 4, $\widehat{\Delta}_M$ is best everywhere, as anticipated, but $\widehat{\Delta}_{mert}$ is close behind, marginally ahead of $\widehat{\Delta}_{gmm}$ even for $n = m = 10,000$, and well ahead for smaller sample sizes. Efficiency comparisons for the convolution of a Normal and an

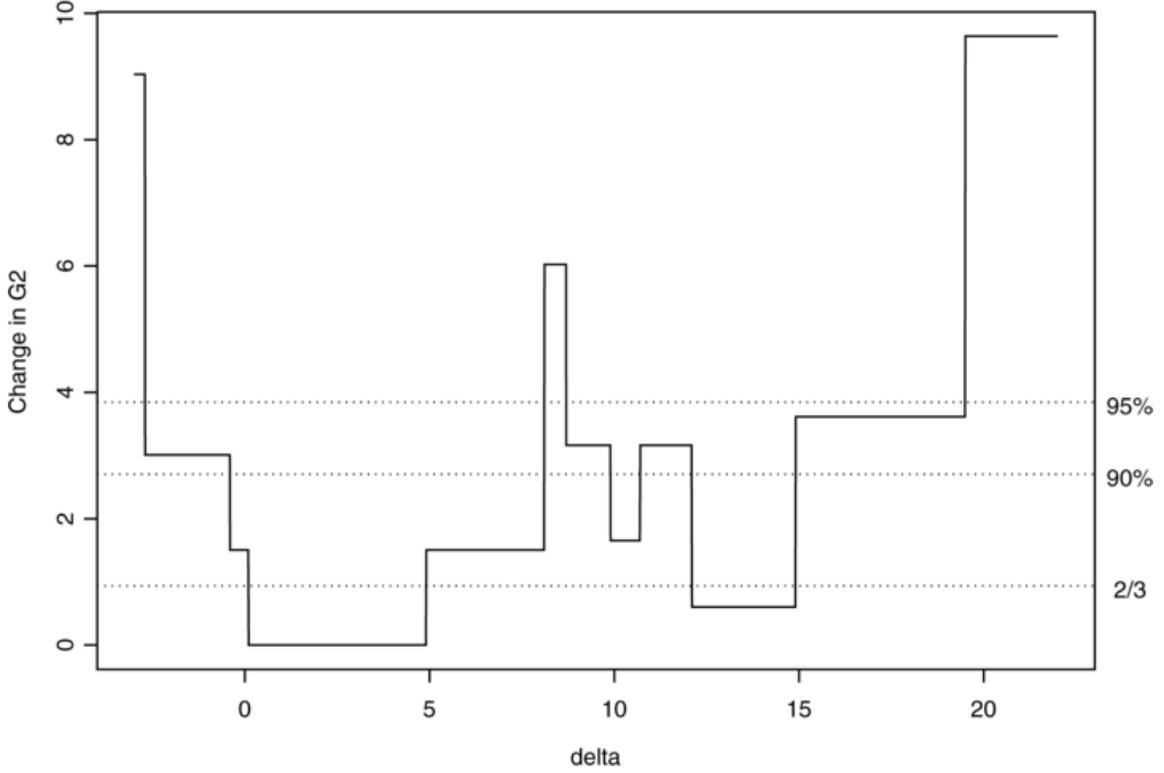

Fig. 4. *Plot of $G^2_{\Delta_0} - G^2_{\widehat{\Delta}_G}$ vs. $\Delta_0$. Note that the set of $\Delta_0$ not rejected (i.e., the confidence set) is not an interval for $\alpha = 0.05$, $0.10$ and $\frac{1}{3}$.*

TABLE 4
*Efficiency for samples from the Cauchy distribution*

| $n$ | **24** | **50** | **20** | **80** | **500** | **2,000** | **10,000** |
|---|---|---|---|---|---|---|---|
| $m$ | **24** | **50** | **80** | **80** | **500** | **2,000** | **10,000** |
| $\widehat{\Delta}_{\mathrm{gmm}} : \widehat{\Delta}_{\mathrm{HL}}$ | 0.82 | 0.92 | 0.90 | 0.95 | 1.07 | 1.11 | 1.18 |
| $\widehat{\Delta}_{\mathrm{gmm}} : \widehat{\Delta}_{\mathrm{M}}$ | 0.66 | 0.75 | 0.74 | 0.76 | 0.86 | 0.90 | 0.94 |
| $\widehat{\Delta}_{\mathrm{gmm}} : \widehat{\Delta}_{\mathrm{mert}}$ | 0.70 | 0.78 | 0.78 | 0.80 | 0.89 | 0.93 | 0.99 |

Summary: $\widehat{\Delta}_{\mathrm{M}}$ is best in all cases, and $\widehat{\Delta}_{\mathrm{mert}}$ is second in all cases.

Exponential random variable are given in Table 5. The best estimator in all cases in Table 5 is $\widehat{\Delta}_{\mathrm{HL}}$. The relative performance of $\widehat{\Delta}_{\mathrm{gmm}}$ improves with increasing sample size, but it is still 5% behind $\widehat{\Delta}_{\mathrm{HL}}$ for $n = m = 10{,}000$. Also, $\widehat{\Delta}_{\mathrm{mert}}$ is ahead of $\widehat{\Delta}_{\mathrm{gmm}}$ up to $n = m = 2{,}000$.

In summary, the relative efficiency of the GMM estimator $\widehat{\Delta}_{\mathrm{gmm}}$ does increase with increasing sample size, as the asymptotic theory says it should, but the improvement is remarkably slow. The fixed score estimator, $\widehat{\Delta}_{\mathrm{mert}}$, is designed to achieve reasonable performance for both the Normal and the Cauchy, and it is more efficient than $\widehat{\Delta}_{\mathrm{gmm}}$ for *all* three sampling distributions in Tables 3 to 5 for sample sizes up to $n = m = 2{,}000$.

## 3.3 Confidence Intervals

In each of the $3 \times 7 = 21$ sampling situations in Tables 3–5, we computed the large sample nomi-

TABLE 5
*Efficiency for samples from the Normal + Exponential distribution*

| $n$ | **24** | **50** | **20** | **80** | **500** | **2,000** | **10,000** |
|---|---|---|---|---|---|---|---|
| $m$ | **24** | **50** | **80** | **80** | **500** | **2,000** | **10,000** |
| $\widehat{\Delta}_{\mathrm{gmm}} : \widehat{\Delta}_{\mathrm{HL}}$ | 0.85 | 0.78 | 0.76 | 0.83 | 0.87 | 0.91 | 0.95 |
| $\widehat{\Delta}_{\mathrm{gmm}} : \widehat{\Delta}_{\mathrm{M}}$ | 1.11 | 1.03 | 0.98 | 1.08 | 1.17 | 1.19 | 1.25 |
| $\widehat{\Delta}_{\mathrm{gmm}} : \widehat{\Delta}_{\mathrm{mert}}$ | 0.94 | 0.85 | 0.83 | 0.91 | 0.95 | 0.99 | 1.04 |

Summary: $\widehat{\Delta}_{\mathrm{HL}}$ is best in all cases, and $\widehat{\Delta}_{\mathrm{mert}}$ is second for $m, n \leq 2000$.

nally 95% confidence intervals from $\widehat{\Delta}_{HL}$, $\widehat{\Delta}_{mert}$ and $\widehat{\Delta}_{gmm}$, and empirically determined the actual coverage rate. For comparison, a binomial proportion with 5000 trials and probability of success 0.95 has standard error 0.003, so $0.95 \pm (2 \times 0.003)$ is 0.944 to 0.956.

The group rank confidence intervals performed well, with coverage close to the nominal level, even when large sample approximations were applied to small samples. All of the $2 \times 21 = 42$ simulated coverage rates for $\widehat{\Delta}_{HL}$ and $\widehat{\Delta}_{mert}$ were between 93.7% and 96.1% and only one was less than 94%. Exact intervals are available for $\widehat{\Delta}_{HL}$ and $\widehat{\Delta}_{mert}$ based on the hypergeometric distribution, but there is no need to simulate these, because their coverage rates are exactly as stated. The empirical coverages of the 95% confidence interval based on GMM are displayed in Table 6. As in Section 3.2, the asymptotic theory appears correct in the limit but takes hold very slowly. The coverage of the nominal 95% interval is about 90% for $n = m = 80$ and about 92% for $n = m = 500$. Also, by the results in Section 3.2, these intervals from $\widehat{\Delta}_{gmm}$ not only have lower coverage than the intervals for the group rank test, but also for sample sizes up to $n = m = 2{,}000$ the intervals from $\widehat{\Delta}_{gmm}$ are typically longer intervals than those from $\widehat{\Delta}_{mert}$ as well. This is not much of a trade-off: lower coverage combined with longer intervals.

Finally, all of the confidence sets from $\widehat{\Delta}_{HL}$ and $\widehat{\Delta}_{mert}$ are intervals, but for $\widehat{\Delta}_{gmm}$ the confidence sets from GMM are often not intervals, and become intervals only by including interior segments that the test rejected. For $\widehat{\Delta}_{gmm}$ for the Normal, with $n = m = 24$, only 49% of the confidence sets are intervals, rising to 67% for $n = m = 500$ and 83% for $n = m = 10{,}000$. Results for the other two distributions are not very different.

### 3.4 GMM's Goodness of Fit Test

Recall that $G^2_{\widehat{\Delta}_{gmm}}$ is often used as a test of the model, in the current context, comparing it to the chi-square distribution on two degrees of freedom. Here, too, the asymptotic properties appear true but are approached very slowly. In particular, when the shift model is correct, $G^2_{\widehat{\Delta}_{gmm}}$ tends to be too small, compared to chi-square on two degrees of freedom, both in the tail and on average. For instance, the chi-square $P$-value from $G^2_{\widehat{\Delta}_{gmm}}$ is less than 0.05 in 0.1% of samples of size $n = m = 24$ from the Normal, in 0.6% of samples of size $n = m = 80$, in 2% of samples of size $n = m = 500$, and 3.7% of samples of size $n = m = 10{,}000$. Similarly, instead of expectation 2 for a chi-square with two degrees of freedom, the expectation of $G^2_{\widehat{\Delta}_{gmm}}$ was 0.9 for samples of size $n = m = 24$ from the Normal, 1.1 for samples of size $n = m = 80$, 1.4 for samples of size $n = m = 500$, and 1.8 for samples of size $n = m = 10{,}000$. Similar results were found for the Cauchy and Normal + Exponential. In sharp contrast, one compares $G^2_{\Delta}$ to chi-square with three degrees of freedom, and $E(G^2_{\Delta}) = 3$ exactly in samples of every size from every distribution; see Section 2.2. In other words, replacing the true $\Delta$ by the estimate $\widehat{\Delta}_{gmm}$ and reducing the degrees of freedom by one to compensate is an adequate correction only in very large samples. To understand the behavior of $G^2_{\widehat{\Delta}_{gmm}}$, it helps to recall what happened in the example in Section 2.4. There, $G^2_{\Delta}$ was minimized at a peculiar choice $\widehat{\Delta}_{gmm}$ of $\Delta$, in part because $G^2_{\Delta}$ avoided not only implausible shifts but also tables like Table 2 which suggest unequal dispersion. Having avoided Table 2—that is, having avoided evidence of unequal dispersion by its choice of $\widehat{\Delta}_{gmm}$—the goodness of fit test, $G^2_{\widehat{\Delta}_{gmm}}$, found no evidence of unequal dispersion. This suggests it may be best to separate two tasks, namely estimation assuming a model is true, and testing the goodness of fit of the model.

## 4. DISPLACEMENT EFFECTS

In the example in Section 1.2, the HL estimate gave a more reasonable estimate of shift than did the GMM estimate assuming the shift model to be true, and based on that estimate, raised clearer doubts about whether the shift model was appropriate. This is seen in Figures 2 and 3. Having raised doubts about the shift model, it is natural to seek exact inferences for the magnitude of the effect without assuming a shift. The shift model is not needed for an exact inference comparing two distributions. There are $112 = 16 \times 7$ possible comparisons of the $n = 16$ exposed subjects to the $m = 7$ controls, and in $V = 87$ of these comparisons the exposed subject had a higher response, where $V$ is the Mann–Whitney statistic. Under the null hypothesis of no treatment effect in a randomized experiment, the chance that $V \geq 82$ is 0.044. It follows from the argument in Rosenbaum (2001, Section 4) that in a randomized experiment, we would be $1 - 0.044 = 95.6\%$ confident that at least $87 - 82 + 1 = 6$ of the 112 possible

TABLE 6
*Empirical coverage of nominal 95% intervals from GMM*

| $n$ | **24** | **50** | **20** | **80** | **500** | **2,000** | **10,000** |
| $m$ | **24** | **50** | **80** | **80** | **500** | **2,000** | **10,000** |
|---|---|---|---|---|---|---|---|
| Normal | 85.7 | 89.7 | 88.0 | 90.1 | 91.9 | 92.3 | 94.2 |
| Cauchy | 87.2 | 90.9 | 89.3 | 91.3 | 92.7 | 93.6 | 94.1 |
| Normal + Exponential | 87.2 | 90.5 | 87.6 | 89.8 | 92.1 | 93.2 | 94.3 |

Summary: 95% intervals from GMM miss too often for $n, m \leq 2,000$.

comparisons, or about 5% of them, favor the exposed group because of effects of the treatment, and the remaining $87 - 6$ favorable comparisons could be due to chance. So the effect is not plausibly zero but could be quite small. The methods in Rosenbaum (2001, Section 5) may be used to display the sensitivity of this inference to departures from random assignment of treatments in an observational study of the sort described in Section 1.2.

## 5. SUMMARY: WHAT ROLE FOR ASSUMPTIONS?

In the radiation effects example in Section 1.2, the generalized method of moments (GMM) estimator, $\widehat{\Delta}_{\mathrm{gmm}}$, estimated a small shift, one that did little to align the boxplots in Figure 2, and the associated test of the shift model using $G^2_{\widehat{\Delta}_{\mathrm{gmm}}}$ suggested that the shift model was plausible. In contrast, the Hodges–Lehmann estimate, $\widehat{\Delta}_{\mathrm{HL}}$, estimated a larger shift, one that did align the centers of the boxplots in Figure 2, but with this shift, the shift model seemed implausible, leading to the analysis in Section 4 which dispensed with the shift model. Although one should not make too much of a single, small example, our sense in this one instance was that GMM's results gave an incorrect impression of what the data had to say.

The simulation considered situations in which the shift model was true. The promise of full asymptotic efficiency with GMM did seem to be true, but was very slow in coming, requiring astonishingly large sample sizes. For moderate sample sizes, $m \leq 500$, $n \leq 500$, GMM was neither efficient nor valid: it produced longer 95% confidence intervals often with coverage well below 95%. In contrast, asymptotic results were a good guide to the performance of the group rank statistics for all sample sizes considered, even for samples of size $n = m = 24$. Moreover, exact inference is straightforward with group rank statistics. The estimator $\widehat{\Delta}_{\mathrm{mert}}$ aims to avoid bad performance under a range of assumptions rather than to achieve optimal performance under one set of assumptions. By every measure, in every situation, $\widehat{\Delta}_{\mathrm{mert}}$ was better than $\widehat{\Delta}_{\mathrm{gmm}}$ for $m \leq 2000$, $n \leq 2000$. Our goal has not been to provide a comprehensive analysis of the radiation effects example, nor to provide improved methods for estimating a shift. Rather, our goal was to create a laboratory environment—transparent, quiet, simple, undisturbed—in which two strategies for creating estimators might be compared. The laboratory conditions were favorable for GMM: (i) the shift parameter is strongly identified, (ii) there are only three moment conditions and (iii) the optimal weight matrix for the moment conditions is known exactly and is free of unknown parameters. In a theorem, the assumptions are the premises of an argument, and for the sole purpose of proving the theorem, they play similar roles: the same conclusion with fewer assumptions is a "better" theorem, or at least better in certain important senses. When used in scientific applications, these same assumptions acquire different roles. As in the quote from Tukey in Section 1.1, assumptions needed for validity of confidence intervals and hypothesis tests play a different role from assumptions used for efficiency or stringency, and both play a very different role from the hypothesis itself. A familiar instance of this arises with hypotheses: omnibus hypotheses (ones that assume very little) are not automatically better hypotheses than focused hypotheses (ones that assume much more)—power may be much higher for the focused hypotheses, and which is relevant depends on the science of the problem at hand. The trade-off discussed by Tukey is a less familiar instance. Here, we have examined a small, clean theoretical case study, in which the same information is used by different methods that embody different attitudes toward assumptions. The group rank statistics followed Tukey's advice, in which validity was obtained by general permutation test arguments, but the weights used in the tests reflected judgements informed by statistical theory in an effort to obtain decent efficiency for a variety of sampling distributions. The generalized method of moments (GMM) tried to estimate the weights, and thereby always have the most efficient procedure, at least asymptotically. In point of fact, the gains in efficiency with GMM did not materialize until very large sample sizes were reached, whereas validity of confidence intervals was severely compromised in samples of conventional size.

## APPENDIX: LARGE SAMPLE THEORY UNDER LOCAL ALTERNATIVES

This appendix discusses the asymptotic efficiency of GMM against local alternatives. Hansen's results about GMM concerned statistics that are differentiable in ways that rank statistics are not, but his conclusions hold nonetheless if differentiability is replaced by asymptotic linearity using Theorem 1 from Jurečková ([1969](#)); see Theorem A.2 below. Here we consider the limiting behavior of $D_{\Delta_N}$ and $G^2_{\Delta_N}$ when $H_0 : \Delta = \Delta_N$ is false but nearly correct, that is, when $F(x) = G(x + \Delta)$ but $\Delta_N = \Delta - \frac{\delta}{\sqrt{N}}$, as $N = n + m \to \infty$ with $\lambda_N = n/(n + m) \to \lambda$, $0 < \lambda < 1$. Because $N \to \infty$, quantities from earlier sections computed from the sample of size $N$ now have an $N$ subscript, for example, $E(\mathbf{A}_{N,\Delta}) = \mathbf{E}_N$ and $\mathrm{var}(\mathbf{A}_{N,\Delta}) = \mathbf{V}_N$. Then $Z^\Delta_{N1}, \ldots, Z^\Delta_{NN}$ are i.i.d. $F(\cdot)$, but $D_{N,\Delta_N}$ and $G^2_{N,\Delta_N}$ are computed from $Z^{\Delta - \delta/\sqrt{N}}_{N1}$, $\ldots, Z^{\Delta - \delta/\sqrt{N}}_{NN}$, so $Z^{\Delta - \delta/\sqrt{N}}_{N1}, \ldots, Z^{\Delta - \delta/\sqrt{N}}_{Nm}$ corresponding to the $X$'s are i.i.d. $F(\cdot)$, but $Z^{\Delta - \delta/\sqrt{N}}_{N,m+1}, \ldots, Z^{\Delta - \delta/\sqrt{N}}_{NN}$ corresponding to the $Y$'s have the same distribution but shifted upwards by $\delta/\sqrt{N}$. Assume $F$ has density $f$ with finite Fisher information and write

$$\varphi(u, f) = -\frac{f'\{F^{-1}(u)\}}{f\{F^{-1}(u)\}}$$

and

$$\eta_g = \lambda(1 - \lambda) \int_{(g-1)/4}^{g/4} \varphi(u, f) \, du \quad \text{for } g = 1, 2, 3, 4.$$

Let $\eta = (\eta_1, \eta_2, \eta_3, \eta_4)^T$ and notice that $0 = \sum \eta_g$ by Hájek, Šidák, and Sen ([1999](#), Lemma 1, page 18). Then a result of Jurečková ([1969](#), Theorem 3.1, page 1891) yields:

THEOREM A.1 (Jurečková, [1969](#)). *For any fixed* $\mathbf{w} = (0, w_2, w_3, w_4)^T$ *with* $0 = w_1 \leq w_2 < w_3 \leq w_4$, $\varepsilon > 0$ *and* $\Upsilon > 0$

$$\lim_{N \to \infty} \Pr\left( \max_{|\delta| < \Upsilon} \left| \frac{1}{\sqrt{N}} \mathbf{w}^T (\mathbf{A}_{N,\Delta - \delta/\sqrt{N}} \right. \right.$$

$$- \mathbf{A}_{N,\Delta}) - \delta \mathbf{w}^T \eta \bigg|$$

$$\left. \geq \varepsilon \sqrt{\frac{1}{N} \mathbf{w}^T \mathbf{V}_N \mathbf{w}} \right) = 0.$$

As $N \to \infty$, one has $(1/N)\mathbf{V}_N \to \mathbf{\Sigma}$ and $N\mathbf{V}^-_N \to \mathbf{\Sigma}^-$ where $\mathbf{\Sigma}$ has entries $\sigma_{jj} = 3\lambda(1 - \lambda)/16$ and $\sigma_{ij} = -\lambda(1 - \lambda)/16$ for $i \neq j$, and $\mathbf{\Sigma}^-$, like $\mathbf{V}^-_N$, has first row and column equal to zero. Moreover, Theorem [A.1](#) implies

$$(9) \quad \frac{\mathbf{w}^T(\mathbf{A}_{N,\Delta - \delta/\sqrt{N}} - \mathbf{E}_N)}{\sqrt{\mathbf{w}^T \mathbf{V}_N \mathbf{w}}} \xrightarrow{D} N\left( \frac{\delta \mathbf{w}^T \eta}{\sqrt{\mathbf{w}^T \mathbf{\Sigma} \mathbf{w}}}, 1 \right).$$

The noncentrality parameter in ([9](#)), $\mathbf{w}^T \eta / \sqrt{\mathbf{w}^T \mathbf{\Sigma} \mathbf{w}}$, is maximized with $\mathbf{w} = \mathbf{\Sigma}^- \eta$, so the best group rank statistic has $\mathbf{w} = \mathbf{\Sigma}^- \eta$. The GMM confidence interval for $\Delta$ was calculated by comparing $G^2_{N,\Delta} - G^2_{N,\widehat{\Delta}_{\mathrm{gmm}}}$ to the chi-square distribution with one degree of freedom. Theorem [A.2](#) shows $G^2_{N,\Delta} - G^2_{N,\widehat{\Delta}}$ converges in probability (henceforth $\xrightarrow{P}$) to the best group rank statistic.

THEOREM A.2. *If $F$ has finite Fisher information, then*

$$(10) \quad \frac{\{(\mathbf{A}_{N,\Delta} - \mathbf{E}_N)^T \mathbf{V}^-_N \eta\}^2}{\eta^T \mathbf{V}^-_N \eta}$$

$$- (G^2_{N,\Delta} - G^2_{N,\widehat{\Delta}_{\mathrm{gmm}}}) \xrightarrow{P} 0.$$

PROOF. Write $\|\mathbf{a}\|_N = \sqrt{N \mathbf{a}^T \mathbf{V}^-_N \mathbf{a}}$. Then $G^2_{N,\Delta} = \|\frac{1}{\sqrt{N}}(\mathbf{A}_{N,\Delta} - \mathbf{E}_N)\|^2_N$ and $G^2_{N,\Delta - \delta/\sqrt{N}} = \|\frac{1}{\sqrt{N}} \cdot (\mathbf{A}_{N,\Delta - \delta/\sqrt{N}} - \mathbf{E}_N)\|^2_N$. Define $\widetilde{G}^2_{N,\Delta - \delta/\sqrt{N}} = \|\frac{1}{\sqrt{N}} \cdot (\mathbf{A}_{N,\Delta} - \mathbf{E}_N) + \delta \eta\|^2_N$, which is quadratic in $\delta$, and is minimized at

$$\widetilde{\delta}_N = \frac{-\sqrt{N}(\mathbf{A}_{N,\Delta} - \mathbf{E}_N)^T \mathbf{V}^-_N \eta}{N \eta^T \mathbf{V}^-_N \eta}$$

$$\doteq \frac{-(1/\sqrt{N})(\mathbf{A}_{N,\Delta} - \mathbf{E}_N)^T \mathbf{\Sigma}^- \eta}{\eta^T \mathbf{\Sigma}^- \eta}$$

$$\xrightarrow{D} N\left( 0, \frac{1}{\eta^T \mathbf{\Sigma}^- \eta} \right);$$

moreover, $\widetilde{G}^2_{N,\Delta-\delta/\sqrt{N}}$ has minimum value $G^2_{N,\Delta} - \{(\mathbf{A}_{N,\Delta} - \mathbf{E}_N)^T \mathbf{V}_N^- \eta\}^2 / \eta^T \mathbf{V}_N^- \eta$. Let $\mathcal{E}_\Upsilon$ be the event $\mathcal{E}_\Upsilon = \{\sqrt{N}|\widehat{\Delta}_{\mathrm{gmm}} - \Delta| < \Upsilon \text{ and } |\widetilde{\delta}_N| < \Upsilon\}$. It is possible to pick a large $\Upsilon > 0$ such that for all sufficiently large $N$, the probability $\Pr(\mathcal{E}_\Upsilon)$ is arbitrarily large. Therefore, in proving (10), we assume $\mathcal{E}_\Upsilon$ has occurred. Write $\psi_{N,\delta} = \{(\mathbf{A}_{N,\Delta-\delta/\sqrt{N}} - \mathbf{A}_{N,\Delta})/\sqrt{N}\} - \delta\eta$ and note that Theorem A.1 implies $\max_{|\delta|<\Upsilon} \|\psi_{N,\delta}\|^2_N \xrightarrow{P} 0$. By the triangle inequality, for any norm $\|\cdot\|$,

$$
(11) \quad \begin{aligned} &\left|\|\mathbf{a}\|^2 - \|\mathbf{b}\|^2\right| \\ &\quad \leq \|\mathbf{a} - \mathbf{b}\|^2 + 2\|\mathbf{a} - \mathbf{b}\| \max(\|\mathbf{a}\|, \|\mathbf{b}\|). \end{aligned}
$$

With $\mathbf{a} = (\mathbf{A}_{N,\Delta-\delta/\sqrt{N}} - \mathbf{E}_N)/\sqrt{N}$ and $\mathbf{b} = (\mathbf{A}_{N,\Delta} - \mathbf{E}_N)/\sqrt{N} + \delta\eta$, so $\mathbf{a} - \mathbf{b} = \psi_{N,\delta}$, then (11) yields

$$
\begin{aligned} &\max_{|\delta|<\Upsilon} |G^2_{N,\Delta-\delta/\sqrt{N}} - \widetilde{G}^2_{N,\Delta-\delta/\sqrt{N}}| \\ &\quad \leq \max_{|\delta|<\Upsilon} \Big\{ \|\psi_{N,\delta}\|^2_N \\ &\qquad + 2\|\psi_{N,\delta}\|_N \\ &\qquad\quad \cdot \max\Big(\sqrt{\widetilde{G}^2_{N,\Delta-\delta/\sqrt{N}}}, \sqrt{G^2_{N,\Delta-\delta/\sqrt{N}}}\Big)\Big\} \\ &\xrightarrow{P} 0, \end{aligned}
$$

so that

$$
\left|\min_{|\delta|<\Upsilon} G^2_{N,\Delta-\delta/\sqrt{N}} - \min_{|\delta|<\Upsilon} \widetilde{G}^2_{N,\Delta-\delta/\sqrt{N}}\right| \xrightarrow{P} 0,
$$

which is (10). $\square$

## ACKNOWLEDGMENTS

Joseph L. Gastwirth and Paul R. Rosenbaum were supported by grants from the National Science Foundation.